%% file: xtalk.tex
\documentclass{article}

\input{diagrams}
\diagramstyle[PostScript=dvips]
\usepackage{epsfig}
\input{xenab}

\input{xpenhdpics}

\title{\fc-multicategories}
\author{Tom Leinster%
	\thanks{Financial support from EPSRC; curveless diagrams by Paul
	Taylor's macros}\\ \\
	\normalsize{Department of Pure Mathematics, University of
	Cambridge}\\ 
	\normalsize{Email: leinster@@dpmms.cam.ac.uk}\\
	\normalsize{Web: http://www.dpmms.cam.ac.uk/$\sim$leinster}}
\date{\normalsize{Notes for a talk at the 70th Peripatetic Seminar on Sheaves
and Logic, Cambridge, 28 February 1999}}

\begin{document}
\maketitle
\begin{abstract}	\centering
What \fc-multicategories are, and two uses for them.
\end{abstract}

\section*{Introduction}
\fc-multicategories are a very general kind of two-dimensional structure,
encompassing bicategories, monoidal categories, double categories and
ordinary multicategories. Here we define what they are and explain how they
provide a natural setting for two familiar categorical ideas. The first is
the \emph{bimodules} construction, traditionally carried out on suitably
cocomplete bicategories but perhaps more naturally carried out on
\fc-multicategories. The second is \emph{enrichment}: there is a theory of
categories enriched in an \fc-multicategory, which includes the usual case
of enrichment in a monoidal category, the obvious extension of this to
ordinary multicategories, and the less well known case of enrichment in
a bicategory.

To finish we briefly indicate the wider context, including how the work below
is just the simplest case of a much larger phenomenon and the reason for the
name `\fc-multicategory'.

\section{What is an \fc-multicategory?}

An \fc-multicategory consists of
\begin{itemize}
\item A collection of \demph{objects} $x$, $x'$, \ldots
\item For each pair \pr{x}{x'} of objects, a collection of \demph{vertical
1-cells} \vslob{x}{}{x'}, denoted $f$, $f'$, \ldots
\item For each pair \pr{x}{x'} of objects, a collection of \demph{horizontal
1-cells} $x\go x'$, denoted $m$, $m'$, \ldots
\item For each $n\geq 0$, objects $\range{x_0}{x_n},x,x'$, vertical 1-cells
$f,f'$, and horizontal 1-cells $\range{m_1}{m_n},m$, a collection of
\demph{2-cells}
\begin{equation}	\label{eq:two-cell}
\begin{diagram}[height=2em]
x_0	&\rTo^{m_1}	&x_1	&\rTo^{m_2}	&\ 	&\cdots	
&\ 	&\rTo^{m_n}	&x_n	\\
\dTo<{f}&		&	&		&\Downarrow	&
&	&		&\dTo>{f'}\\
x	&		&	&		&\rTo_{m}	&	
&	&		&x',	\\
\end{diagram}
\end{equation}
denoted $\theta$, $\theta'$, \ldots
\item \demph{Composition} and \demph{identity} functions making the objects
and vertical 1-cells into a category
\item A \demph{composition} function for 2-cells, as in the picture
\begin{diagram}[width=.5em,height=1em]
\blob&\rTo^{m_1^1}&\cdots&\rTo^{m_1^{r_1}}&
\blob&\rTo^{m_2^1}&\cdots&\rTo^{m_2^{r_2}}&\blob&
\ &\cdots&\ &
\blob&\rTo^{m_n^1}&\cdots&\rTo^{m_n^{r_n}}&\blob\\
\dTo<{f_0}&&\Downarrow\theta_1&&
\dTo&&\Downarrow\theta_2&&\dTo&
\ &\cdots&\ &
\dTo&&\Downarrow\theta_n&&\dTo>{f_n}\\
\blob&&\rTo_{m_1}&&
\blob&&\rTo_{m_2}&&\blob&
\ &\cdots&\ &
\blob&&\rTo_{m_n}&&\blob\\
\dTo<{f}&&&&&&&&\Downarrow\theta &&&&&&&&\dTo>{f'}\\
\blob&&&&&&&&\rTo_{m}&&&&&&&&\blob\\
\end{diagram}
$\goesto$
\begin{diagram}[width=.5em,height=1em]
\blob&\rTo^{m_1^1}&\ &&
&&&&\cdots&
&&&
&&\ &\rTo^{m_n^{r_n}}&\blob\\
&&&&&&&&&&&&&&&&\\
\dTo<{f\of f_0}&&&&&&&&\Downarrow\theta\of\tuple{\theta_1}{\theta_2}{\theta_n}
&&&&&&&&\dTo>{f'\of f_n}\\
&&&&&&&&&&&&&&&&\\
\blob&&&&&&&&\rTo_{m}&&&&&&&&\blob\\
\end{diagram}
($n\geq 0, r_i\geq 0$, with \blob's representing objects) 
\item An \demph{identity} function
\[
\begin{diagram}[width=1em,height=2em]
x&\rTo^{m}&x'\\
\end{diagram}
\diagspace\goesto\diagspace
\begin{diagram}[width=1em,height=2em]
x&\rTo^{m}&x'\\
\dTo<{1_x}&\Downarrow 1_{m}&\dTo>{1_{x'}}\\
x&\rTo_{m}&x'\\
\end{diagram}
\]
\end{itemize}
such that 2-cell composition and identities obey associativity and identity
laws. 

\begin{eg}
\item
Any double category gives an \fc-multicategory, in which a 2-cell as at
\bref{eq:two-cell} is a 2-cell
\begin{diagram}[height=2em]
x_0	&\rTo^{m_n \of \cdots \of m_1}	&x_n\\
\dTo<{f}&\Downarrow			&\dTo>{f'}\\
x	&\rTo_{m}			&x'\\
\end{diagram}
in the double category.

\item
Any bicategory gives an \fc-multicategory in which the only vertical 1-cells
are identity maps, and a 2-cell as at \bref{eq:two-cell} is a 2-cell
\[
x_0
\ctwo{m_n \of\cdots\of m_1}{m}{}
x_n
\]
in the bicategory (with $x_0=x$ and $x_n=x'$).

\item	\label{eg:mon-cat}
Any monoidal category gives an \fc-multicategory in which there is one
object and one vertical 1-cell, and a 2-cell
\begin{equation}	\label{eq:vt-two-cell}
\begin{diagram}[size=2em,abut]
\bullet	&\rLine^{M_1}	&\bullet	&\rLine^{M_2}	&\bullet	&\cdots 
&\bullet	&\rLine^{M_n}	&\bullet	\\
\dLine<1&		&	&		&\Downarrow&
&	&		&\dLine>1\\
\bullet	&		&	&		&\rLine_{M}&
&	&		&\bullet\\
\end{diagram}
\end{equation}
is a morphism $M_n\otimes\cdots\otimes M_1 \go M$.

\item
Similarly, any ordinary multicategory gives an \fc-multicategory: there is
one object, one vertical 1-cell, and a 2-cell \bref{eq:vt-two-cell} is a map
$\range{M_1}{M_n}\go M$.

\item	\label{span}
We define an \fc-multicategory \Span. Objects are sets, vertical 1-cells are
functions, a horizontal 1-cell $X\go Y$ is a diagram
$
\begin{diagram}[width=1.8em,height=1em]
	&	&M	&	&	\\
	&\ldTo	&	&\rdTo	&	\\
X	&	&	&	&Y	\\
\end{diagram},
$
and a 2-cell inside
\begin{equation}	\label{eq:span-two-cell}
\begin{diagram}[width=1.8em,height=1em,tight]
	&	&M_1	&	&	&	&M_2	&	&	
&	&	&	&M_n	&	&	\\
	&\ldTo	&	&\rdTo	&	&\ldTo	&	&\rdTo	&
&\cdots	&	&\ldTo	&	&\rdTo	&	\\ 
X_0	&	&	&	&X_1	&	&	&	&\ 
&	&\ 	&	&	&	&X_n	\\
	&	&	&	&	&	&	&	&
&	&	&	&	&	&	\\
\dTo<f	&	&	&	&	&	&	&	&
&	&	&	&	&	&\dTo>{f'}\\
	&	&	&	&	&	&	&M	&
&	&	&	&	&	&	\\
	&	&	&	&	&	&\ldTo(7,2)&	&\rdTo(7,2)
&	&	&	&	&	&	\\
X	&	&	&	&	&	&	&	&
&	&	&	&	&	&X'	\\
\end{diagram}
\end{equation}
is a function $\theta$ making
\begin{diagram}[width=2em,height=1em]
	&	&M_n\of\cdots\of M_1	&	&	\\
	&\ldTo	&			&\rdTo	&	\\
X_0	&	&			&	&X_n	\\
	&	&\dTo>{\theta}		&	&	\\
\dTo<f	&	&			&	&\dTo>{f'}\\
	&	&M			&	&	\\
	&\ldTo	&			&\rdTo	&	\\
X	&	&			&	&X'	\\
\end{diagram}
commute, where $M_n\of\cdots\of M_1$ is the limit of the top row of
\bref{eq:span-two-cell}. Composition is defined in the obvious way.

\end{eg}

\section{Bimodules}

Bimodules have traditionally been discussed in the context of
bicategories. Thus given a bicategory \cat{B}, we construct a new bicategory
\Bim{\cat{B}} whose 1-cells are bimodules in \cat{B} (see \cite{CKW} or
\cite{Kos}). The drawback is that to do this, we must make certain assumptions
about the behaviour of local coequalizers in \cat{B}.

However, the \fcat{Bim} construction extends to \fc-multicategories, and
working in this context allows us to drop all the technical assumptions: we
therefore obtain a functor $\fcat{Bim}: \fc\hyph\Multicat \go
\fc\hyph\Multicat$. The definition is rather dry, so we omit it here and just
give a few examples; the reader is referred to \cite[2.6]{GECM} for further
details. 

\begin{eg}

\item	\label{eg:ab}
Let $V$ be the \fc-multicategory coming from the monoidal category
\pr{\Ab}{\otimes} (see \bref{eg:mon-cat} above). Then \Bim{V} has
\begin{description}
\item[objects:]
rings
\item[vertical 1-cells:]
ring homomorphisms
\item[horizontal 1-cells $R\go S$:]
\pr{S}{R}-bimodules
\item[2-cells:]
A 2-cell
\begin{diagram}[height=2em]
R_0	&\rTo^{M_1}	&R_1	&\rTo^{M_2}	&\ 	&\cdots	
&\ 	&\rTo^{M_n}	&R_n	\\
\dTo<{f}&		&	&		&\Downarrow\,\theta&
&	&		&\dTo>{f'}\\
R	&		&	&		&\rTo_{M}	&	
&	&		&R'	\\
\end{diagram}
is a multi-additive map $M_n\times\cdots\times M_1 \goby{\theta} M$ of abelian
groups such that
\begin{eqnarray*}
\theta(r_n.m_n, m_{n-1}, \ldots)	&=&
f(r_n).\theta(m_n, m_{n-1}, \ldots)	\\
\theta(m_n.r_{n-1}, m_{n-1}, \ldots)	&=&
\theta(m_n, r_{n-1}.m_{n-1}, \ldots)
\end{eqnarray*}
etc.
\end{description}

\item
If $V$ is the \fc-multicategory \Span\ then \Bim{V} has
\begin{description}
\item[objects:]
monads in \Span, i.e.\ small categories
\item[vertical 1-cells:]
functors
\item[horizontal 1-cells $\scat{C} \go \scat{C'}$:]
profunctors (i.e.\ functors $\scat{C}^{\op} \times \scat{C'} \go \Set$)
\item[2-cells:]
A 2-cell
\begin{diagram}[height=2em]
\scat{C}_0	&\rTo^{M_1}	&\scat{C}_1	&\rTo^{M_2}	&\ 	&\cdots	
&\ 	&\rTo^{M_n}	&\scat{C}_n	\\
\dTo<{F}&		&	&		&\Downarrow	&
&	&		&\dTo>{F'}\\
\scat{C}	&		&	&		&\rTo_{M}	&	
&	&		&\scat{C}'	\\
\end{diagram}
is a natural family of functions
\[
M_1(c_0,c_1) \times\cdots\times M_n(c_{n-1},c_n) \go M(Fc_0,F'c_n),
\]
one for each $c_0 \in \scat{C}_0$, \ldots, $c_n \in \scat{C}_n$.
\end{description}

\item
Let $V$ be the \fc-multicategory coming from a bicategory \Bee\ with
nicely-behaved local coequalizers. If we discard the non-identity vertical
1-cells from \Bim{V} then we obtain the \fc-multicategory coming from the
traditional bicategory \Bim{\Bee}---e.g.\ in \bref{eg:ab}, we get the
bicategory of rings and bimodules. 

\end{eg}

\section{Enrichment}

We define what a `category enriched in $V$' is, for any \fc-multicategory
$V$. This generalizes the established definitions for monoidal categories and
bicategories.

Fix an \fc-multicategory $V$. A \demph{category $C$ enriched in $V$} consists
of
\begin{itemize}
\item 
a set $C_0$ (`of objects')
\item 
for each $a\in C_0$, an object $\eend{C}{a}$ of $V$
\item 
for each $a,b\in C_0$, a horizontal 1-cell $\eend{C}{a}
\rTo^{\ehom{C}{a}{b}} \eend{C}{b}$ in $V$
\item 
for each $a,b,c\in C_0$, a `composition' 2-cell
\begin{diagram}[height=2em]
\eend{C}{a}&\rTo^{\ehom{C}{a}{b}}&\eend{C}{b}	&\rTo^{\ehom{C}{b}{c}}
&\eend{C}{c}	\\
\dTo<{1}&		&\Downarrow\,\comp_{a,b,c}	&	&\dTo>{1}\\
\eend{C}{a}&		&\rTo_{\ehom{C}{a}{c}}		&	&\eend{C}{c}\\
\end{diagram}
\item
for each $a\in C_0$, an `identity' 2-cell
\begin{diagram}[height=2em]
\eend{C}{a}	&\rEquals		&\eend{C}{a}	\\
\dTo<{1}	&\Downarrow\,\id_{a}	&\dTo>{1}	\\
\eend{C}{a}	&\rTo_{\ehom{C}{a}{a}}	&\eend{C}{a}	\\
\end{diagram}
(where the equality sign along the top denotes a string of 0 horizontal
1-cells)
\end{itemize}
such that \comp\ and \id\ satisfy associativity and identity axioms.

\emph{Remark:}\/ We haven't used the vertical 1-cells of $V$ in any significant
way, but we would do if we went on to talk about \emph{functors} between
enriched categories (which we won't here).

\begin{eg}

\item
Let $V$ be (the \fc-multicategory coming from) a monoidal category. Then the
choice of $\eend{C}{a}$'s is uniquely determined, so we just have to specify
the set $C_0$, the \ehom{C}{a}{b}'s, and the maps $\ehom{C}{b}{c} \otimes
\ehom{C}{a}{b} \go \ehom{C}{a}{c}$ and $I \go \ehom{C}{a}{a}$. This gives the
usual notion of enriched category.

\item
If $V$ is an ordinary multicategory then we obtain an obvious generalization
of the notion for monoidal categories: so a category enriched in $V$ consists
of a set $C_0$, an object \ehom{C}{a}{b} of $V$ for each $a,b$, and suitable
maps $\ehom{C}{a}{b}, \ehom{C}{b}{c} \go \ehom{C}{a}{c}$ and $\cdot \go
\ehom{C}{a}{a}$ (where $\cdot$ denotes the empty sequence).

\item
If $V$ is a bicategory then we get the notion of Walters et al
(see \cite{BCSW}, \cite{CKW}, \cite{Wal}).

\item
Let $D$ be a category enriched in \pr{\Ab}{\otimes}. Then we get a category
$C$ enriched in \Bim{\Ab}:
\begin{itemize}
\item $C_0 = D_0$ ($=$ objects of $D$)
\item \eend{C}{a} is the ring \ehom{D}{a}{a} (whose multiplication is
composition in $D$)
\item \ehom{C}{a}{b} is the abelian group \ehom{D}{a}{b} acted on by
$\eend{C}{a} = \ehom{D}{a}{a}$ (on the right) and $\eend{C}{b} =
\ehom{D}{b}{b}$ (on the left)
\item composition and identities are as in $D$.
\end{itemize}
So the passage from $D$ to $C$ is basically down to the fact that composition
makes \ehom{D}{a}{a} into a ring and \ehom{D}{a}{b} into a
\pr{\ehom{D}{b}{b}}{\ehom{D}{a}{a}}-bimodule. It's a very mechanical process,
and in fact for general $V$ there's a functor
\[
(\mbox{categories enriched in }V) \go (\mbox{categories enriched 
in }\Bim{V}). 
\]

\item
An example of a category enriched in \Bim{\Span} ($=$ categories $+$ functors
$+$ profunctors\ldots): $C_0$ is \nat, \eend{C}{n} is the category of
$n$-dimensional real differentiable manifolds and diffeomorphisms, and the
profunctor \ehom{C}{m}{n} is the functor
\begin{eqnarray*}
\eend{C}{m}^{\op} \times \eend{C}{n}	&\go&
\Set	\\
\pr{M}{N}				&\goesto&
\{\mbox{differentiable maps }M\go N\}.
\end{eqnarray*}

\item
Let \fcat{ParBjn} be the sub-\fc-multicategory of \Span\ in which all
horizontal 1-cells are of the form $(X \lMonic M \rMonic Y)$: so this 1-cell
is a partial bijection between $X$ and $Y$. Let $S$ be a set and $(C_i)_{i\in
I}$ a family of subsets. Then we get a category $C$ enriched in
\fcat{ParBjn}:
\begin{itemize}
\item $C_0=I$
\item $\eend{C}{i} = C_i$
\item $\ehom{C}{i}{j} = (C_i \lMonic C_i \cap C_j \rMonic C_j)$
\item $\comp_{i,j,k}$ is the inclusion $C_i \cap C_j \cap C_k \sub
C_i \cap C_k$
\item $\id_{i}$ is the inclusion $C_i \sub C_i \cap C_i$.  
\end{itemize}

\item
Fix a topological space $A$. Suppose $A$ is nonempty and path-connected;
choose a basepoint $a_0$ and a path $\gamma_a: a_0 \go a$ for each $a\in
A$. Then we get a category $C$ enriched in the homotopy bicategory $V$ of $A$
(where $V$ consists of points of $A$, paths in $A$, and homotopy classes of
path homotopies in $A$):
\begin{itemize}
\item $C_0=A$
\item $\eend{C}{a}=a$
\item $\ehom{C}{a}{b}$ is
$\begin{array}{c}\epsfig{file=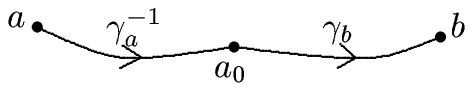}\end{array}$
\item composition $\ehom{C}{b}{c}\of\ehom{C}{a}{b} \go \ehom{C}{a}{c}$ is the
(homotopy class of the) obvious homotopy from
\[
\begin{array}{c}\epsfig{file=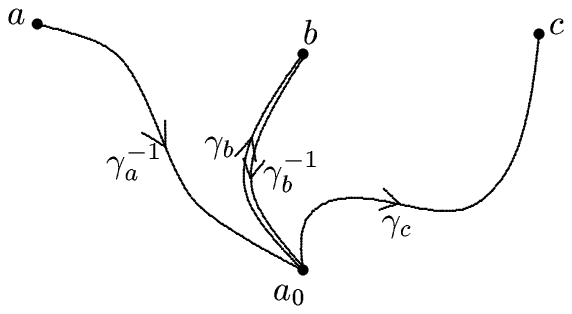}\end{array}
\mbox{to}
\begin{array}{c}\epsfig{file=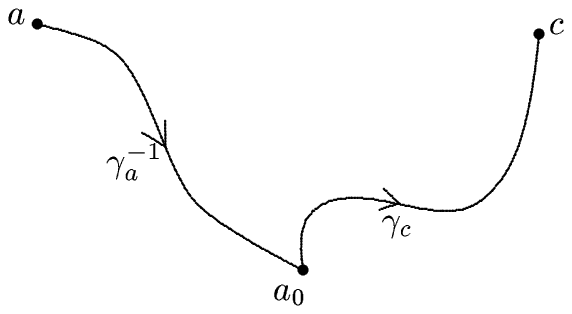}\end{array}
\]
\item identities work similarly.
\end{itemize}
\end{eg}

\section*{The wider context}

Given a monad $T$ on a category \Eee\ (with certain properties), there's a
category of \demph{$T$-multicategories} (see \cite{GOM}, \cite{Bur} or
\cite{Her}). For example:
\begin{quote}
\begin{tabular}{c|c}
\Cartpr				&$T$-multicategories	\\
\hline
\Zeropr				&categories				\\
\pr{\Set}{\mbox{free monoid}}	&ordinary multicategories		\\
\pr{\Gph}{\mbox{\textbf{\underline{f}}ree \textbf{\underline{c}}ategory}}&
\fc-multicategories	
\end{tabular}
\end{quote}
where $\Gph=\ftrcat{(\bullet\pile{\rTo\\ \rTo}\bullet)}{\Set}$.

Moreover, if one defines a \demph{$T$-graph} to be a diagram 
$
\begin{diagram}[width=1.5em,height=1em]
	&	&C_1	&	&	\\
	&\ldTo	&	&\rdTo	&	\\
T(C_0)	&	&	&	&C_0	\\
\end{diagram}
$
in \Eee, then there's a forgetful functor
$
T\hyph\Multicat \go T\hyph\Gph,
$
this has a left adjoint, and the adjunction is monadic. Write $\Eeep =
T\hyph\Gph$ and $T'$ for the induced monad on \Eeep. Then we can also
discuss $T'$-multicategories, and in fact there's a theory of
$T$-multicategories enriched in a $T'$-multicategory.

The simplest case is $\Cartpr = \Zeropr$: then $\Cartprp =
\pr{\Gph}{\mr{free\ category}}$, so we have a theory of categories enriched
in an \fc-multicategory. This is just the theory we discussed above. 

A full explanation of these ideas can be found in \cite{GECM}.

\end{document}

%% file: xenab.tex
\newcommand{\mcm}[3]{\newcommand{#1}[#2]{{\ensuremath{#3}}}}
\usepackage{latexsym}
\usepackage{amstex}
\mcm{\hyph}{0}{\mbox{-}}
\mcm{\diagspace}{0}{\mbox{\hspace{2em}}}
\mcm{\mc}{1}{\mathcal{#1}}
\mcm{\mr}{1}{\mathrm{#1}}
\mcm{\mi}{1}{\mathit{#1}}
\mcm{\mb}{1}{\mathbf{#1}}
\mcm{\cat}{1}{\mc{#1}}
\mcm{\scat}{1}{\mathbb{#1}}
\mcm{\fcat}{1}{\mb{#1}}
\newcommand{\url}[1]{\mbox{\tt #1}}
\mcm{\blob}{0}{\raisebox{0.2mm}{\ensuremath{\scriptscriptstyle\bullet}}}
\mcm{\of}{0}{\raisebox{0.2mm}{\ensuremath{\scriptstyle\circ}}}
\mcm{\op}{0}{\mr{op}}
\mcm{\id}{0}{\mi{id}}
\mcm{\Set}{0}{\fcat{Set}}
\mcm{\Multicat}{0}{\fcat{Multicat}}
\mcm{\Graph}{0}{\fcat{Graph}}
\mcm{\nat}{0}{\mathbb{N}}	
\mcm{\pr}{2}{\tuplebts{#1,#2}}
\mcm{\range}{2}{#1,\,\ldots\,,#2}
\mcm{\tuplebts}{1}{(#1)}
\mcm{\tuple}{3}{\tuplebts{\range{#1,#2}{#3}}}
\mcm{\bftuple}{2}{\tuplebts{\range{#1}{#2}}}
\mcm{\ftrcat}{2}{[#1,#2]}
\mcm{\go}{0}{\rTo}
\mcm{\goby}{1}{\rTo^{#1}}
\mcm{\goesto}{0}{\,\longmapsto\,}
\mcm{\vslob}{3}
	{\left.
	\begin{diagram}[height=1.5em]
	#1		\\
	\dTo>{\,#2}	\\
	#3		\\
	\end{diagram}
	\right.}
\newarrow{Equals}=====
\newarrow{Monic}{vee}---{vee}
\newarrow{Mod}--+->
\mcm{\sub}{0}{\,\subseteq\,}
\mcm{\ob}{1}{\mr{ob}\,#1}
\mcm{\comp}{0}{\mi{comp}}
\mcm{\Eee}{0}{\cat{E}}
\mcm{\Eeep}{0}{\cat{E'}}
\mcm{\Cartpr}{0}{\pr{\Eee}{T}}
\mcm{\Cartprp}{0}{\pr{\Eeep}{T'}}
\mcm{\Gph}{0}{\fcat{Graph}}
\mcm{\fc}{0}{\fcat{fc}}
\mcm{\Zeropr}{0}{\pr{\Set}{\id}}
\mcm{\Bim}{1}{\fcat{Bim}(#1)}
\mcm{\Bee}{0}{\cat{B}}
\mcm{\Ab}{0}{\fcat{Ab}}
\mcm{\Span}{0}{\fcat{Span}}
\newcommand{\bref}[1]{(\ref{#1})}
\mcm{\ehom}{3}{#1[#2,#3]}
\mcm{\eend}{2}{#1[#2]}

\newcommand{\demph}[1]{\textbf{#1}}
\newenvironment{eg}{\paragraph{Examples}\begin{enumerate}}{\end{enumerate}}
\newarrow{Inc}{boldhook}--->

%% file: xpenhdpics.tex
\newlength{\gwidth}	
\newlength{\gvert}	
\newlength{\gdrop}	
\newlength{\gbaredrop}	
\newlength{\goffset}	
\newlength{\gtemp}	
\newcommand{\present}[1]{%
\makebox[1\gwidth]{%
\rule[-1\gdrop]{0ex}{1\gvert}%
\raisebox{-1\gbaredrop}{#1}}}
\newcommand{\cinitdims}[2]{%
\setlength{\unitlength}{1em}%
\setlength{\goffset}{.35\unitlength}%
\setlength{\gwidth}{#1\unitlength}%
\setlength{\gvert}{#2\unitlength}%
\setlength{\gdrop}{.5\gvert}%
\addtolength{\gdrop}{-1\goffset}%
\setlength{\gbaredrop}{1\gdrop}%
\addtolength{\gvert}{.6\unitlength}%
\addtolength{\gdrop}{.3\unitlength}}	
\newcommand{\abovepic}[1]{%
\settoheight{\gtemp}{\ensuremath{#1}}%
\addtolength{\gvert}{1\gtemp}%
\settodepth{\gtemp}{\ensuremath{#1}}%
\addtolength{\gvert}{1\gtemp}}
\newcommand{\belowpic}[1]{%
\settoheight{\gtemp}{\ensuremath{#1}}%
\addtolength{\gvert}{1\gtemp}%
\addtolength{\gdrop}{1\gtemp}%
\settodepth{\gtemp}{\ensuremath{#1}}%
\addtolength{\gvert}{1\gtemp}%
\addtolength{\gdrop}{1\gtemp}}
\newcommand{\cell}[4]{\put(#1,#2){\makebox(0,0)[#3]{\ensuremath{#4}}}}
\newcommand{\prectwo}[3]%
{\begin{picture}(4.2,3.4)(-0.1,-0.2)%
\cell{2}{3.2}{b}{#1}%
\cell{2}{-0.2}{t}{#2}%
\cell{2.2}{1.5}{l}{#3}%
\qbezier(0,2)(2,4)(4,2)%
\qbezier(0,1)(2,-1)(4,1)%
\put(4,2){\vector(1,-1){0}}%
\put(4,1){\vector(1,1){0}}%
\cell{2}{1.5}{c}{\Downarrow}%
\end{picture}}
\mcm{\ctwo}{3}{%
\cinitdims{4.2}{3.4}%
\abovepic{#1}%
\belowpic{#2}%
\present{\prectwo{#1}{#2}{#3}}}